\documentclass[12pt]{article}
\newcommand{\ol}{\setlength{\itemsep}{0pt.}\begin{enumerate}}
\newcommand{\eol}{\end{enumerate}\setlength{\itemsep}{-\parsep}}

\newtheorem{thm}{Theorem}

\newtheorem{pr}{Proposition}
\newtheorem{co}{Corollary}

\usepackage{amsmath, amsfonts}

\newcommand{\e}{\varepsilon}
\newcommand{\N}{\mathbb{N}}

\def\X{{\frak X}}

\def\pf{ \noindent {\bf Proof: \  }}

\newcommand{\qed}{\hfill\vrule height6pt
width6pt depth0pt}

\def\endpf{\qed \medskip}

\renewcommand{\qed}{\hfill\vrule height6pt  width6pt depth0pt}

\def\ss#1{_{\lower3pt\hbox{$\scriptstyle #1$}}}

\title{{The number of closed ideals in $L(L_p)$  }
\author{William B. Johnson\thanks{Supported in part by NSF DMS-1900612}\  \
Gideon Schechtman\thanks{Supported in part by Israel Science
Foundation}
}}
\begin{document}
\maketitle

\begin{abstract}
We show that there are $2^{2^{\aleph_0}}$ different closed ideals in the Banach algebra $L(L_p(0,1))$, $1<p\not= 2<\infty$.
This solves a problem in A. Pietsch's 1978 book ``Operator Ideals". The proof is quite different from other methods of producing closed ideals in the space of bounded operators on a Banach space; in particular, the ideals are not contained in the strictly singular operators and yet do not contain projections onto subspaces that are non Hilbertian. We give a criterion for a space with an unconditional basis to  have $2^{2^{\aleph_0}}$ closed ideals in terms of the existence of a single operator on the space with some special asymptotic properties. We then show that for $1<q<2$ the space ${\frak X}_q$ of Rosenthal, which is isomorphic to a complemented subspace of $L_q(0,1)$, admits such an operator.
\end{abstract}

\section{Introduction}
For a reasonably complete discussion of the history of constructing closed ideals in $L(L_p)$, see the introduction in \cite{jps}. Here we just remark that in 1981, Bourgain, Rosenthal, and the second author \cite{brs} constructed $\aleph_1$ mutually non isomorphic complemented subspaces of $L_p:= L_p(0,1)$ for $1<p\not= 2 < \infty$, thereby producing (as noted in \cite{pietsch}) $\aleph_1$ different closed ideals in $L(L_p)$. (It is of course well known that the compact operators are the only closed ideal in $L(L_2)$.) At that time it was open whether, absent the continuum hypothesis, $L(L_p)$ contains a continuum of closed ideals.  Recently, Schlumprecht and Zs\'ak \cite{sz} built a continuum of closed ideals in $L_p:= L_p(0,1)$.

The main contribution of this paper is Theorem \ref{thm:main} in which we prove that
$L(L_p)$, $1<p\not= 2 < \infty$, has exactly $2^{2^{\aleph_0}}$ different closed ideals.

Recall the notions of small and large closed ideal in $L(X)$. An ideal is called small if it is contained in the ideal of strictly singular operators. Otherwise it is called large. The ideals built in \cite{sz} are all small, while the ones coming from infinite dimensional  complemented subspaces are clearly large. Our basic construction is designed to produce large ideals. Note that there are at most a continuum of non mutually isomorphic complemented subspaces of $L_p$ (as the density character of $L(L_p)$ and of the set of projections on $L_p$ is the continuum). So necessarily we produce  different  kinds of ideals. Unfortunately, we   do not produce any new complemented subspaces of $L_p$.

The new large ideals in $L(L_p)$ that we construct are ``smallish" in the sense that, even though there are idempotents in the ideals whose ranges are isomorphic to    $\ell_2$ (see Remark 5), no operator in any of the ideals is an isomorphism on a copy of $\ell_p$. The Kadec--Pe\l czy\'nski dichotomy principle \cite{kp}   implies   that every complemented subspace of $L_p$ that is not isomorphic to a Hilbert space contains a complemented subspace that is isomorphic to $\ell_p$. Consequently, the range of any infinite rank idempotent in any of the ideals that we construct in Theorem \ref{thm:main} (and, as we said, there are infinite rank idempotents in the ideals) {\sl must} be isomorphic to $\ell_2$.

To put these new ``smallish" large ideals into perspective within the Banach algebra $L(L_p)$, notice that it follows from the Kadec--Pe\l czy\'nski dichotomy principle \cite{kp} that  there are exactly two different minimal large closed ideals in $L(L_p)$ when $2<p < \infty$, and thus also for $1<p<2$ (because an operator $T$  in $L(L_p)$ is strictly singular if and only if $T^*$ is strictly singular on $L_q$, $1/p+1/q=1$, by Weis' theorem \cite{wei}).  The first of these is $\Gamma_{\ell_p}(L_p)$, the ideal of operators that factor through $\ell_p$. This ideal is closed because an operator  $T:X\to L_p$, $2<p<\infty$, factors through $\ell_p$ if and only if  $I_{p,2}T$  is compact, where $I_{p,2}$ is the formal identity mapping from $L_p$ into $L_2$; see \cite{joh}. One can prove using the Kadec--Pe\l czy\'nski dichotomy principle \cite{kp} that $I_{p,2}S$ is compact whenever $S$ is a strictly singular operator on $L_p$, so the alternate characterization of $\Gamma_{\ell_p}(L_p)$ for $2<p<\infty$ also yields that  $\Gamma_{\ell_p}(L_p)$ contains all strictly singular operators on $L_p$, $2<p<\infty$,  and thus also for $1<p<2$ by \cite{wei}.

The second minimal large closed ideal in $L(L_p)$ is the closure $\overline{\Gamma_2}(L_p)$
 of the ideal $\Gamma_2(L_p)$ of operators on $L_p$ that factor through a Hilbert space. Here the closure is needed; in fact, it is not hard to see that there are compact operators on $L_p$ that do not factor through a Hilbert space.

 We recall in passing that as was noted in \cite{jps} the situation in $L(L_1)$ is nicer: $\Gamma_{\ell_1}(L_1)$ is the unique minimal closed large ideal in $L(L_1)$ and it contains all the strictly singular operators on $L_1$.

 In Remark 6 we prove that the new large ideals we construct in $L(L_p)$ do not contain the strictly singular operators on $L_p$, and hence neither does $\overline{\Gamma_2}(L_p)$. All previously known large ideals in $L(L_p)$ other than $\overline{\Gamma_2}(L_p)$ do contain the strictly singular operators, and this is the first proof that $\overline{\Gamma_2}(L_p)$ does not. A byproduct of Remark 6, stated as Remark 7, is that $L(L_p)$ contains exactly $2^{2^{\aleph_0}}$ small closed ideals.

Our construction and proof of Theorem \ref{thm:main} consist of two steps. In Section \ref{sec:proposition} we state and prove the technical Proposition \ref{pr:main}. This easily yields Corollary \ref{co:main},  which gives a general criterion for a space with an unconditional basis to contain $2^{2^{\aleph_0}}$ different closed ideals. The criterion is in term of the existence of a special operator on the space.

In Section \ref{sec:the operator} we show that for $1<q<2$, the space  $L_q$ contains a complemented subspace (this is Rosenthal's $\X_q$ space, which has an unconditional basis) that admits an operator satisfying the criterion of Proposition \ref{pr:main}. The construction here borrows a lot from a previous similar construction from \cite{jsmult}. Duality and complementation then imply the main result.

\section{The main proposition}\label{sec:proposition}

There is a continuum   of infinite subsets  of the natural numbers $\mathbb{N}$ each two of which have only finite intersection. Denote some fixed such continuum by $\mathcal{C}$. For a finite dimensional normed space $E$, we denote by $d(E)$ the Banach--Mazur distance (isomorphism constant) of $E$ to a Euclidean space. Also, recall that, for an operator $T:X\to Y$ between two normed spaces, $\gamma_2(T)$ denotes its factorization constant through a Hilbert space:
\[
\gamma_2(T)=\inf\{\|A\|\|B\|\ ;\  T=AB, A:H\to Y, B:X\to H,\ H\  \mbox {a Hilbert space}\}.
\]
If $T$ is of rank $k$, then $\gamma_2(T)\le k^{1/2}\|T\|$ because every $k$ dimensional normed space is $k^{1/2}$-isomorphic to $\ell_2^k$ \cite[Theorem 15.5]{T-J}. Note that $d(E)$ is just $\gamma_2(I_E)$, where $I_E$ is the identity operator on $E$.

\begin{pr}\label{pr:main}
Let $X$ be a Banach space with a $1$-unconditional basis $\{e_i\}$, let $Y$ be a Banach space, and let $T:X\to Y$ be an operator of norm at most one satisfying:\\ \\
(a) For some $\eta>0$ and for every $M$ there is a finite dimensional subspace $E$ of $X$ such that $d(E)>M$ and $\|Tx\|> \eta\|x\|$ for all $x\in E$.\\ \\
(b) For some constant $\Gamma$ and every $m$ there is an $n$ such that every $m$-dimensional subspace $E$ of $[e_i]_{i\ge n}$ satisfies $\gamma_2(T_{|E})\le \Gamma$.\\ \\
Then there exist natural numbers $1=p_1<q_1<p_2<q_2<\dots$ such that, denoting for each $k$, $G_k:=[e_i]_{i=p_k}^{q_k}$, and  defining for each $\alpha\in \mathcal{C}$, the operator $P_\alpha:X\to[G_k]_{k\in\alpha}$ to be the the natural basis projection, and setting $T_\alpha:=TP_\alpha$, we have the following:\\ \\
If $\alpha_1,\dots,\alpha_s\in \mathcal{C}$ (possibly with repetitions) and $\alpha\in \mathcal{C}\setminus\{\alpha_1,\dots,\alpha_s\}$, then for all $A_1,\dots,A_s\in L(Y)$ and all $B_1,\dots,B_s\in L(X)$,
\begin{equation}\label{eq:pr1}
\|T_\alpha-\sum_{i=1}^s A_iT_{\alpha_i}B_i\|\ge \eta/2.
\end{equation}
\end{pr}

\pf Note first that we can strengthen condition (a) to include also that given any $n$ one can chose the subspace $E$ to also satisfy that it is contained in $[e_i]_{i>n}$.
Now choose inductively $1=p_1<q_1<p_2<q_2\dots$ so that for each $k$, $G_k=[e_i]_{i=p_k}^{q_k}$ contains a subspace $E_k$ with $\|Tx\|> \eta\|x\|$ for all $x\in E_k$ and
\[
d(E_k)\ge q_{k-1}
\]
(as we'll see, it is enough that $d(E_k)/q_{k-1}^{1/2}\to \infty$)
and, if $E$ is a subspace of $H_k=[G_l]_{l=p_{k+1}}^\infty$ with $\dim E\le q_k$, then
\[
\gamma_2(T_{|E})<\Gamma.
\]
Let now $P_\alpha:X\to[G_k]_{k\in\alpha}$ be the natural basis projection and set $T_\alpha:=TP_\alpha$.\\
Suppose that  $\alpha_1,\dots,\alpha_s\in \mathcal{C}$ (possibly with repetitions) and $\alpha\in \mathcal{C}\setminus\{\alpha_1,\dots,\alpha_s\}$. Assume to the contrary that there are $A_1,\dots,A_s\in L(Y)$ and  $B_1,\dots,B_s\in L(X)$ such that
\begin{equation}\label{eq:approx}
\|T_\alpha-\sum_{i=1}^s A_iT_{\alpha_i}B_i\|< \eta/2.
\end{equation}
There are infinitely many $k\in \alpha\setminus\bigcup_{i=1}^s\alpha_i$. For each such $k$ let $R_k$ be the basis projection onto $[G_l]_{l<k}$ and $Q_k$ the basis projection onto $[G_l]_{l>k}$. Now for any $i=1,\dots,s$ we have $T_{\alpha_i}G_k=0$ since $k\notin \alpha_i$,  and  $\dim(R_kB_iE_k)\le q_{k-1}$ and $\dim(B_iE_k)\le q_k$, so we get that for each $i$,
\begin{eqnarray*}
\gamma_2(A_i T_{\alpha_i}B_{i|E_k})
&\le& \gamma_2(A_i T_{\alpha_i}R_k B_{i|E_k})+\gamma_2(A_i T_{\alpha_i}Q_k B_{i|E_k})\\
&\le& q_{k-1}^{1/2}\|A_i\|\|B_i\|+\Gamma\|A_i\|\|B_i\|.
\end{eqnarray*}
Consequently,
\begin{equation}\label{eq:upperbound}
\gamma_2(\sum_{i=1}^s A_i T_{\alpha_i}B_{i|E_k})\le (\max_{1\le i\le s}\|A_i\|\|B_i\|)s(q_{k-1}^{1/2}+\Gamma).
\end{equation}
On the other hand, since $\|x\|
\ge \|T_\alpha x\|\ge \eta\|x\|$ for all $x\in E_k$, (\ref{eq:approx}) implies that
\[
(1+\eta/2)\|x\|\ge\|\sum_{i=1}^s A_i T_{\alpha_i}B_i x\|\ge\eta\|x\|/2
\]
for all $x\in E_k$. Since $d(E_k)\ge q_{k-1}$,  we deduce that
\[
\gamma_2(\sum_{i=1}^s A_i T_{\alpha_i}B_{i|E_k})\ge \frac{\eta}{2+\eta}q_{k-1}.
\]
For $k$ large enough this contradicts (\ref{eq:upperbound}).
\endpf

\noindent{\bf Remark 0.} Observe that the only condition on $T_\alpha$ that was used to get the inequality (\ref{eq:pr1}) is that $\|x\| \ge \|T_\alpha x\| \ge \eta\|x\|$ for all $x$ in $E_k$ with $k\in \alpha$. Consequently, the proof of Corollary \ref{co:main} below shows that any  operator $S$ in $L(X)$ for which there is $\eta > 0$ such that $  \|S x\| \ge \eta\|x\|$ for all $x$ in $E_k$ with $k\in \alpha$ cannot be  in the closed  ideal generated  by $\{T_\beta: \beta \in \mathcal{C}, \beta\not= \alpha\}$.  In fact, from the proof of Proposition \ref{pr:main}, only the inequality $\|S x\| \ge \eta\|x\|$ for all $x$ in $H_k$ with $k\in \alpha$ and where  $H_k$ is  isomorphic to   $E_k$ with isomorphism constant independent of $k$ is sufficient to conclude that $S$ is not  in the closed  ideal generated  by $\{T_\beta: \beta \in \mathcal{C}, \beta\not= \alpha\}$.
This observation is used in Remark 6 at the end of this paper.

\begin{co}\label{co:main}
Let $X$ be a Banach space with a $1$-unconditional basis $\{e_i\}$ and assume there is an operator $T:X\to X$ of norm at most one satisfying $(a)$ and $(b)$ of Proposition \ref{pr:main}. Then $L(X)$ has exactly $2^{2^{\aleph_0}}$ different closed ideals.
\end{co}

\pf
For any nonempty proper subset $\mathcal{A}$ of $\mathcal{C}$ let $\mathcal{I}_{\mathcal{A}}$ be the ideal generated by $\{T_\alpha\}_{\alpha\in\mathcal{A}}$; i.e., all operators of the form $\sum_{i=1}^s A_iT_{\alpha_i}B_i$ with $s\in\mathbb{N}$, $A_i,B_i\in L(X)$, $\alpha_i\in \mathcal{A}$, $i=1,\dots,s$.  To avoid cumbersome  notation, interpret
$\mathcal{A}\subset \mathcal{C}$ to mean that $\mathcal{A}$ is a nonempty proper subset of
$\mathcal{C}$.

Since we allow repetition of the $T_{\alpha_i}$, it is easy to see that this really defines a (non closed) ideal. Let $\mathcal{B}$ be a subset of $\mathcal{C}$ different from  $\mathcal{A}$ and assume, without loss of generality, that $\mathcal{B}\not\subset\mathcal{A}$. Let $\alpha\in \mathcal{B}\setminus\mathcal{A}$. Then by  Proposition \ref{pr:main},  $T_{\alpha}\notin \overline{\mathcal{I}_{\mathcal{A}}}$. Consequently, $\{\overline{\mathcal{I}_{\mathcal{A}}}\}_{\mathcal{A}\subset\mathcal{C}}$ are all different.

Since the density character of $L(X)$, for any separable $X$, is at most the continuum, it is easy to see that, for any separable space $X$, $L(X)$ has at most $2^{2^{\aleph_0}}$ different closed ideals.
\endpf


\noindent{\bf Remarks:}

\noindent {\bf 1.} One can strengthen the conclusion of the corollary by getting an antichain of $2^{2^{\aleph_0}}$ closed ideals in $L(X)$; i.e., such a collection no two of whose members are included one in the other. For that one just uses a collection of $2^{2^{\aleph_0}}$ subsets of $\mathcal{C}$ no two of which are included one in the other.

\noindent{\bf 2.} Similarly, one gets a collection of $2^{\aleph_0}$ different closed ideals in $L(X)$ that form a chain (by taking a chain of subsets of $\mathcal{C}$ of that cardinality). It is also easy to show by a density argument that, for any separable $X$, this is the maximal cardinality of any chain of closed ideals in $L(X)$.

\noindent{\bf 3.} If $Y$ is a Banach space that contains a complemented subspace $X$ with the properties of Corollary \ref{co:main} then clearly $L(Y)$ also has  $2^{2^{\aleph_0}}$ different closed ideals (actually an antichain). The same is true also for any space isomorphic to such a $Y$.

\noindent{\bf 4.} The simplest examples of spaces $X$ that satisfy the hypotheses of Corollary \ref{co:main} and thus $L(X)$ has $2^{2^{\aleph_0}}$ different closed ideals are $(\sum \ell_{r_i}^{n_i})_2$ for $r_i \uparrow 2$ and $n_i$ satisfying $n_i^{\frac{1}{r_i}-\frac12}\to\infty$. Consequently,  by  Remark 3,  $L((\sum \ell_{r_i})_2)$ for $r_i \uparrow 2$ also has $2^{2^{\aleph_0}}$ different closed ideals.
Interesting, but less natural, examples of separable spaces $X$ with $L(X)$ having $2^{2^{\aleph_0}}$ different closed ideals were known before (see \cite{man}). Unfortunately
$(\sum \ell_{r_i}^{n_i})_2$  for $r_i\uparrow 2$ and $n_i^{\frac{1}{r_i}-\frac12}\to\infty$ does not embed isomorphically as a complemented subspace into any $L_p$, $p<\infty$, so this example is not good for our purposes. Actually, at least for some sequences $\{(r_i,n_i)\}$ with the above properties, $(\sum \ell_{r_i}^{n_i})_2$ does not even embed isomorphically into any $L_p$ space, $p<\infty$. That this is true, for example,  if each $(r,n)\in \{(r_i,n_i)\}$ repeats $n$ times  follows from Corollary 3.4 in \cite{ks}.

In the next section we show how to get complemented subspaces of the reflexive $L_p$ spaces that satisfy the hypotheses of Corollary \ref{co:main}.

\section{The Operator $T$}\label{sec:the operator}

In this section we prove that for each $1<q<2$ there is a complemented subspace of $L_q$ isomorphic to a space $X$ with a $1$-unconditional basis on which there is an operator of norm at most one with properties $(a)$ and $(b)$ of Proposition \ref{pr:main}.

Recall that for a sequence $u=\{u_j\}_{j=1}^\infty$ of positive real numbers and for $p>2$, the Banach space $\X_{p,u}$ is the  sequence space with norm
\begin{equation}\label{eq:Xp}
\|\{a_j\}_{j=1}^\infty\|=\max\{(\sum_{j=1}^\infty|a_j|^p)^{1/p},(\sum_{j=1}^\infty |a_j u_j|^2)^{1/2}\}.
\end{equation}
Rosenthal \cite{rosxp} proved that $\X_{p,u}$ is isomorphic to a complemented subspace of $L_p$ with  the isomorphism constant and the complementation constant depending only on $p$. If $u$ is such that  for all $\varepsilon>0$, $\sum_{u_j<\varepsilon}u_j^{\frac{2p}{p-2}}=\infty$, then one gets a space isomorphically different from $\ell_p, \ell_2$ and $\ell_p\oplus\ell_2$. However, for different $u$ satisfying the  condition above, the different $\X_{p,u}$ spaces are mutually isomorphic. We denote by $\X_p$ any of these spaces. Later we shall need more properties of the spaces $\X_{p,u}$ and of particular embeddings of them into  $L_p$, but for  now we only need the representation (\ref{eq:Xp}) and we think of $\X_{p,u}$ as a subspace of $\ell_p\oplus_\infty\ell_2$.

Let $\{e_j\}_{j=1}^\infty$ be the unit vector basis of $\ell_p$ and let   $\{f_j\}_{j=1}^\infty$ be the unit vector basis of $\ell_2$. Let $v=\{v_j\}_{j=1}^\infty$ and $w=\{w_j\}_{j=1}^\infty$ be two positive real sequences such that $\delta_j=w_j/v_j\to 0$ as $j\to\infty$ and $\max_{1\le j<\infty}\delta_j\le 1$. Set
\[
g_j^v=e_j+v_jf_j\in\ell_p\oplus_\infty\ell_2  \ \ \mbox{and}\ \ g_j^w=e_j+w_jf_j\in\ell_p\oplus_\infty\ell_2.
\]
Then $\{g_j^v\}_{j=1}^\infty$ is the unit vector basis of $\X_{p,v}$ and $\{g_j^w\}_{j=1}^\infty$ is the unit vector basis of $\X_{p,w}$. Define also
\[
\Delta:\X_{p,w}\to \X_{p,v}
\]
by
\[
\Delta g_j^w=\delta_j g_j^v.
\]
Note that $\Delta$ is the restriction to $\X_{p,w}$ of
$
K\in L(\ell_p\oplus_\infty\ell_2)
$
defined by
\[
K(e_j)=\delta_je_j\ \
\mbox{and}\ \ K(f_j)=f_j
\]
Consequently, $\|\Delta\|\le\|K\|= 1$.

The following proposition follows immediately from the easily verified  fact that \break $\|K_{|[e_j]_{j=m}^\infty}\|\to 0$ as $m\to\infty$.

\begin{pr}\label{pr:Delta} Given $n$ there exists an $m$ such that if $E$ is an $n$ dimensional subspace of ${[e_j]_{j=m}^\infty\oplus[f_j]_{j=1}^\infty} \subset \ell_p\oplus_\infty\ell_2$,
 then $\gamma_2(K_{|E})\le 2$. In particular, if $E$ is an $n$ dimensional subspace of
 ${[g^w_j]_{j=m}^\infty} \subset \X_{p,w}$, then $\gamma_2(\Delta_{|E})\le 2$.
\end{pr}

Next we  define weights $\{v_j\}$ and $\{w_j\}$ with some additional properties. For that we use     different representations of the spaces $\X_{p,u}$. It was proved in \cite{rosxp} that if $\{X_j\}_{j=1}^\infty$, is a sequence of symmetric, each three valued, independent random variables all $L_p$ normalized,
 $2<p<\infty$, then $\{X_j\}_{j=1}^\infty$ is equivalent, in $L_p$, to $\{g_j^u\}_{j=1}^\infty$, the unit vector basis of $\X_{p,u}$, where $u_j=\|X_j\|_2$. Defining $Y_j=X_j/\|X_j\|_q$, for $q=p/(p-1)$, $\{Y_j\}_{j=1}^\infty$ is equivalent, in $L_q$, to the basis $\{h_j^u\}_{j=1}^\infty$ of $\X_{q,u}:=\X_{p,u}^*$ that is dual to the unit vector basis of $\X_{p,u}$.

 Let us say already at this early stage that, for some appropriate weights $\{v_j\}$ and $\{w_j\}$, the operator $T$ we are after will be of the form $\Delta^*$ followed by a norm one isomorphism from $\X_{q,w}$ to $\X_{q,v}$.

Recall that $P:L_p\to [X_j]_{j=1}^\infty$ defined by
\[
Pf=\sum_{j=1}^\infty (\int_0^1 fY_j)X_j
\]
defines a bounded projection onto $[X_j]_{j=1}^\infty$ (and $P^*$ a bounded projection from $L_q$ onto $[Y_j]_{j=1}^\infty$). The norms of the equivalences above and of the projections depend on $p$ but not on the particular weights $u$.

We now recall a construction from Section 4 of \cite{jsmult}. It was shown there that, given $1<q<2$,
any sequence $\{\delta_i\}_{i=1}^\infty$ that decreases to zero, any sequence  $\{r_i\}_{i=1}^\infty$
such that $q<r_i\uparrow 2$ fast enough and in particular satisfying
$\delta_i^{\frac{q(2-r_i)}{2-q}}>1/2$, $i=1,2,\dots$, and for any sequence
$\e_i\downarrow 0$, we can find two sequences $\{Y_i\}$ and
$\{Z_i\}$  of symmetric,
independent, three valued random variables, all normalized in $L_q$, with the following additional
properties:

\begin{itemize}
\item Put $v_j=1/\|Y_j\|_2$ and $w_j=1/\|Z_j\|_2$.
Then there are disjoint finite subsets $\sigma_i$,
$i=1,2,\dots$, of the integers such that $w_j=\delta_i v_j$ for $j\in \sigma_i$.
\item There are independent random variables $\{\bar Y_i\}$ and $\{\bar Z_i\}$,
with $\bar Y_i$  normalized in $L_q$ and $r_i$--stable;  $\bar Z_i$ is $r_i$--stable with $1\ge \|\bar Z_i\|_q\ge 3/4$ for each $i$, and there are coefficients
$\{a_j\}$ such that
\begin{equation}\label{eq:stables}
\|\bar Y_i-\sum_{j\in\sigma_i}a_jY_j\|_q<\e_i \ \ \ \mbox {and}\ \
\ \|\bar Z_i-\sum_{j\in\sigma_i}\delta_ia_jZ_j\|_q<\e_i.
\end{equation}
\end{itemize}

We may of course repeat each of the triplets $r_i,\delta_i,\e_i$-s as many (finitely many)  times as we wish. Thus we conclude that given any sequence $\{\delta_i\}_{i=1}^\infty$ decreasing to zero, any sequence  $\{r_i\}_{i=1}^\infty$
such that $q<r_i\uparrow 2$ and satisfying
$\delta_i^{\frac{q(2-r_i)}{2-q}}>1/2$, $i=1,2,\dots$,  any sequence of integers $n_i$, and any sequence
$\e_i\downarrow 0$, we can find two sequences $\{Y_i\}$ and
$\{Z_i\}$ of symmetric,
independent, three valued random variables, all normalized in $L_q$, with the following additional
properties:
\begin{itemize}
\item Put $v_j=1/\|Y_j\|_2$ and $w_j=1/\|Z_j\|_2$.
Then there are disjoint finite subsets $\sigma_{i,l}$,
$i=1,2,\dots, l=1,\dots n_i$  of the integers such that $w_j=\delta_i v_j$ for $j\in \sigma_{i,l}$.
\item There are independent random variables $\{\bar Y_{i,l}\}$ $r_i$--stable normalized in $L_q$, $\{\bar Z_{i,l}\}$ $r_i$--stable with $1\ge \|\bar Z_{i,l}\|_q\ge 3/4$ for each $i$ and $l$,
 and there are coefficients
$\{a_j\}$ such that
\begin{equation}\label{eq:stables2}
\|\bar Y_{i,l}-\sum_{j\in\sigma_{i,l}}a_jY_j\|_q<\e_i \ \ \ \mbox {and}\ \
\ \|\bar Z_{i,l}-\sum_{j\in\sigma_{i,l}}\delta_ia_jZ_j\|_q<\e_i.
\end{equation}
\end{itemize}

Choosing the $\e_i$ small enough, we can assume that $\{\sum_{j\in\sigma_{i,l}}a_jY_j\}_{l=1}^{n_i}$ is,  in $L_q$,  $2$-equivalent to the unit vector basis of $\ell_{r_i}^{n_i}$, and similarly $\{\sum_{j\in\sigma_{i,l}}\delta_i a_jZ_j\}_{l=1}^{n_i}$ is,  in $L_q$,  $2$-equivalent to the unit vector basis of $\ell_{r_i}^{n_i}$. Denoting by $R$ the map that sends $Y_j$ to $\delta_i Z_j$ for $j\in\sigma_{i,l}$, we get that this map satisfies that for all $i$ there is a space $E_i$ that is  $2$-isomorphic to  $\ell_{r_i}^{n_i}$ such that $\|Rx\|\ge\|x\|/4$ for all $x\in E_i$. Choosing the $n_i$ large enough, we can also assume that for all $k$,
\[
n_i^{\frac{1}{r_i}-{\frac12}} \to \infty  \ \ \mbox{as} \ \ i\to\infty.
\]
Since $n_i^{\frac{1}{r_i}-{\frac12}}$ is the distance of $\ell_{r_i}^{n_i}$ to a Hilbert space, we get that $d(E_i)\to\infty$.

We are now ready to state and prove the main proposition of this section.

\begin{pr}\label{pr:T}
With the choice of $v=\{v_j\}$ and $w=\{w_j\}$ above,  set $X=\X_{q,v}$, let $\Delta^*:\X_{q,v}\to \X_{q,w}$ be the adjoint of $\Delta$ defined at the beginning of this section, and let  $S$ be a norm one isomorphism from $\X_{q,w}$ onto $\X_{q,v}$. Put $T=S\Delta^*$. Then $X,T$ satisfy the assumptions of Proposition \ref{pr:main}.
\end{pr}

\pf
Since $T=ARB$ for isomorphisms $A$ and $B$, the discussion above provides a proof of property $(a)$. Property $(b)$ follows by duality from Proposition \ref{pr:Delta}. Indeed,
fix $m$ and $n$ and let $E$ be an $m$-dimensional subspace of $[h_i^v]_{i\ge n}$. $\Delta^*(E)$ is a subspace of $[h_i^w]_{i\ge n}$, so there is a $k=k(m)$-dimensional subspace $F$ of $[g_i^w]_{i\ge n}$ that $2$-norms $E$. Here, $k=k(m)$ depends only on $m$ (and we used the $1$-unconditionality of the bases). By Proposition \ref{pr:Delta}, for some $n$ depending only on $k$ and thus only on $m$, $\gamma_2(\Delta_{|F})\le 2$. From this it is easy to get that $\gamma_2(\Delta^*_{|E})\le 4$. Consequently, this holds also for $T=S\Delta^*$.
\endpf

\section{The main result and additional comments}

\begin{thm}\label{thm:main}
For every $1<p\not=2<\infty$ the number of different closed ideals in $L(\X_p)$ and in $L(L_p)$ is exactly $2^{2^{\aleph_0}}$. Moreover, each of these spaces contains an antichain of closed ideals of cardinality $2^{2^{\aleph_0}}$ and a chain of cardinality $2^{\aleph_0}$.
\end{thm}

\pf
For $\X_q$, $1<q<2$, the theorem follows from Proposition \ref{pr:T} and Corollary \ref{co:main}. For $\X_p$, $2<p<\infty$, it follows by simple duality. Since for $1<p\not=2<\infty$ the space $\X_p$ is isomorphic to a complemented subspace of $L_p$,  it follows also for $L_p$.

The statements about chains and antichains follow from the remarks at the end of Section \ref{sec:proposition}.
\endpf

\noindent{\bf Remark 5.} As is stated in the introduction, the new ideals in $L(L_p)$ and $ L(\X_p)$, $1<p\not=2<\infty$,  constructed  in Theorem \ref{thm:main} are all large and in fact contain projections whose ranges are isomorphic to $\ell_2$.
\medskip

\pf
First we observe that it is enough to show that for each $\alpha\in \mathcal{C}$, the operator  $T_\alpha$ on $X$ (recall that $X$ is isomorphic to  $\X_q$, where $1<q<2$),  isomorphically preserves a copy of $\ell_2$. Here $T$ is the operator produced in Proposition \ref{pr:T} and $T_\alpha$ is defined in the statement of Proposition \ref{pr:main}. Indeed, since any subspace of $L_q$, $1<q<2$, that is isomorphic to $\ell_2$ contains a further infinite dimensional subspace that  is complemented in $L_q$ (this  fact was probably first observed
by Pe\l czy\'nski; see  \cite[p. 1106]{jsmult} for a proof), this will show that the identity on $\ell_2$ factors through $T_\alpha$  and hence there is a projection in the ideal generated by $T_\alpha$   whose range is isomorphic to $\ell_2$. This will give Remark 5 for
$L(\X_p)$ when $1<p<2$ and the case of $L(\X_p)$ for $2<p<\infty$  follows by duality.
The statement for $L(L_p)$, $1<p\not=2<\infty$, is then immediate.

To show that $T_\alpha$ isomorphically preserves a copy of $\ell_2$, note that
the space $\X_{q,v}$ we built contains a modular space \cite[Def. 4.d.1]{lt1}  $\ell_{\{r_i\}}$ with $r_i\uparrow 2$ on which $T_\alpha$ is an isomorphism and thus (by passing to a subsequence of the sequence $r_i$ that tends quickly to $2$), also contains  an isomorph of $\ell_2$ on which $T_\alpha$ is an isomorphism.
\endpf

\noindent{\bf Remark 6.}
The large ideals in $L(L_q) $ and $L(\X_q)$  constructed in Theorem \ref{thm:main}  do not contain the ideal of strictly singular operators.

\pf (sketch):
By \cite{wei} and how we constructed the ideals in  $L(L_q) $ from the ideals in $L(\X_q)$,  it is enough to consider the ideals constructed in $L(\X_q)$ for $1<q<2$.
Let $T$ be the operator  and $X$ be the space  isomorphic to $\X_q$ that are defined in Proposition \ref{pr:T} and which satisfy the assumptions of Proposition \ref{pr:main}. Let
$\{T_\alpha: \alpha \in  \mathcal{C}\}$ be the corresponding operators on $X$ given by Proposition \ref{pr:T}. As in the proof of Corollary \ref{co:main}, for $ \mathcal{A}$ a (always nonempty, proper) subset of $ \mathcal{C}$ let
$ \mathcal{I}_\mathcal{A}$   be the ideal in $L(X)$ generated by $\mathcal{A}$.  Given $ \mathcal{A} \subset   \mathcal{C}$, take any $\alpha \in  \mathcal{C}$ that is not in $ \mathcal{A}$. We know that $T_\alpha$ is not in $ \overline{\mathcal{I}}_\mathcal{A}$, but   we want a strictly singular operator that is not in $ \overline{\mathcal{I}}_\mathcal{A}$ and $T_\alpha$ is not strictly singular. Let $Y:=( \sum_{k=1}^\infty  G_k)_q$, where the $G_k$ are the block subspaces of $X$  defined in the proof of Proposition \ref{pr:main}. The $G_k$ are contractively complemented in $X$ and $X$ is isomorphic to a complemented subspace of $L_q$, hence $Y$ is isomorphic to a complemented subspace of $\ell_q$ (and thus to $\ell_q$ by Pe\l czy\'nski's well-known theorem, but we do not need this) which in turn is isomorphic to a complemented subspace of $X$. Define $U:Y \to X$ by   making $U$ the identity on each
$G_k$ and extending by linearity and continuity.  This is OK because $L_q$ has type $q$ and
$(G_k)$  is a monotonely unconditional Schauder decomposition for a subspace of $X$, hence the decomposition $(G_k)$ has an upper $q$-estimate (even with constant $1$). Let $E_k$ be the subspace of $G_k$ defined in the proof of Proposition \ref{pr:main}. The operator $T_\alpha U$ is   strictly singular and  $\|T_\alpha U x\| \ge \eta\|x\|$ for all $x$ in $E_k$ with $k$ in $\alpha$. Since $Y$ is isomorphic to a complemented subspace of $X$, we also get a strictly singular operator $S:X\to X$ and subspaces $H_k$ of $Y$ with $H_k$ isomorphic to $E_k$ (with isomorphism constant independent of $k$) such $\|Sx\| \ge \eta \|x\|$ for all $x\in H_k$ with $k\in \alpha$. By Remark 0 after Proposition \ref{pr:main}, this is enough to yield that $S$ is not in the closed ideal in $L(X)$ generated by $\{T_\beta: \beta \in \mathcal{C}, \beta \not=\alpha\}$.
\endpf

\noindent{\bf Remark 7.} $L(L_q) $ and $L(\X_q)$, $1<q \not=2<\infty$ both  contain exactly  $2^{2^{\aleph_0}}$
closed small ideals.

\pf (sketch): Again, it is enough to deal with the case of $L(\X_q)$ with $1<q<2$. Let $X$ and $T$ be as in Remark 6. For $ \mathcal{A} \subset  \mathcal{C}$, let $ \mathcal{J}_ \mathcal{A}$ be the ideal
in $L(X)$ generated by $\{T_\alpha UP: \alpha \in  \mathcal{A} \}$, where $P$ is any fixed projection from $X$ onto a subspace isomorphic to $Y$ (we identify $Y$ with that subspace).
All $ \overline{\mathcal{J}}_\mathcal{A}$ are small ideals and clearly  $ \mathcal{J}_\mathcal{A}$ is contained in the   ideal ${I}_\mathcal{A}$ generated by $\{T_\alpha : \alpha \in  \mathcal{A} \}$. But in Remark 6 we saw that $T_\alpha UP$ is not contained in $\overline{I}_\mathcal{A}$ when $\alpha\not\in \mathcal{A}$, so  $ \overline{\mathcal{J}}_\mathcal{A}\not= \overline{\mathcal{J}}_\mathcal{B}$ when $A\not= B$.
\endpf

\bigskip

\noindent William B. Johnson\newline
             Department Mathematics\newline
             Texas A\&M University\newline
             College Station, TX, USA\newline
             E-mail: johnson@math.tamu.edu

\bigskip

\noindent Gideon Schechtman\newline Department of
Mathematics\newline Weizmann Institute of Science\newline Rehovot,
Israel\newline E-mail: gideon@weizmann.ac.il

\end{document}